%% file: paper.arxiv.tex
\documentclass{acmsiggraph}                     
\pdfoutput=1

\usepackage[scaled=.92]{helvet}
\usepackage{times}

\usepackage{graphicx}

\usepackage{parskip}

\usepackage[labelfont=bf,textfont=it]{caption}

\usepackage{epstopdf}
\usepackage{mathptmx}
\usepackage{epsfig}
\usepackage{graphicx}
\usepackage{amssymb}
\usepackage{multirow}
\usepackage{algorithm}
\usepackage{algorithmic}
\usepackage{subfigure}
\usepackage{color}
\usepackage{amsmath}
\usepackage{wrapfig}
\usepackage{psfrag}

\newif\ifsgpversion
 \sgpversionfalse 
\newif\ifsupplemental
 \supplementalfalse 
\newif\iffullimages
 \fullimagestrue  
\newif\ifconvergencediscussion
 \convergencediscussiontrue  
\newif\ifjakeversion
 \jakeversiontrue  
\newif\ifcleanratio
 \cleanratiotrue  

\newcommand{\ignore}[1]{}

\newcommand{\measure}[1]{d\mu_#1}
\newcommand{\R}[0]{{\mathbb R}}
\newcommand{\tinyZ}{{\hbox{\tiny 0}}}

\newcommand{\dg}[0]{g_\tinyZ^{\hbox{\tiny -1}}g_t}
\newcommand{\dgi}[0]{g_t^{\hbox{\tiny-1}}g_\tinyZ}
\newcommand{\dhg}[0]{h^{\hbox{\tiny -1}}g_t}
\newcommand{\dgh}[0]{g_t^{\hbox{\tiny-1}}h}

\onlineid{papers\_0217}

\title{Can Mean-Curvature Flow Be Made Non-Singular?}

\author{Michael Kazhdan\thanks{e-mail: misha@cs.jhu.edu}\qquad Jake Solomon\thanks{e-mail: jake@math.huji.ac.il}\qquad Mirela Ben-Chen\thanks{e-mail: mirelab@gmail.com}}


\begin{document}

\teaser
{
	\includegraphics[width=1.95\columnwidth]{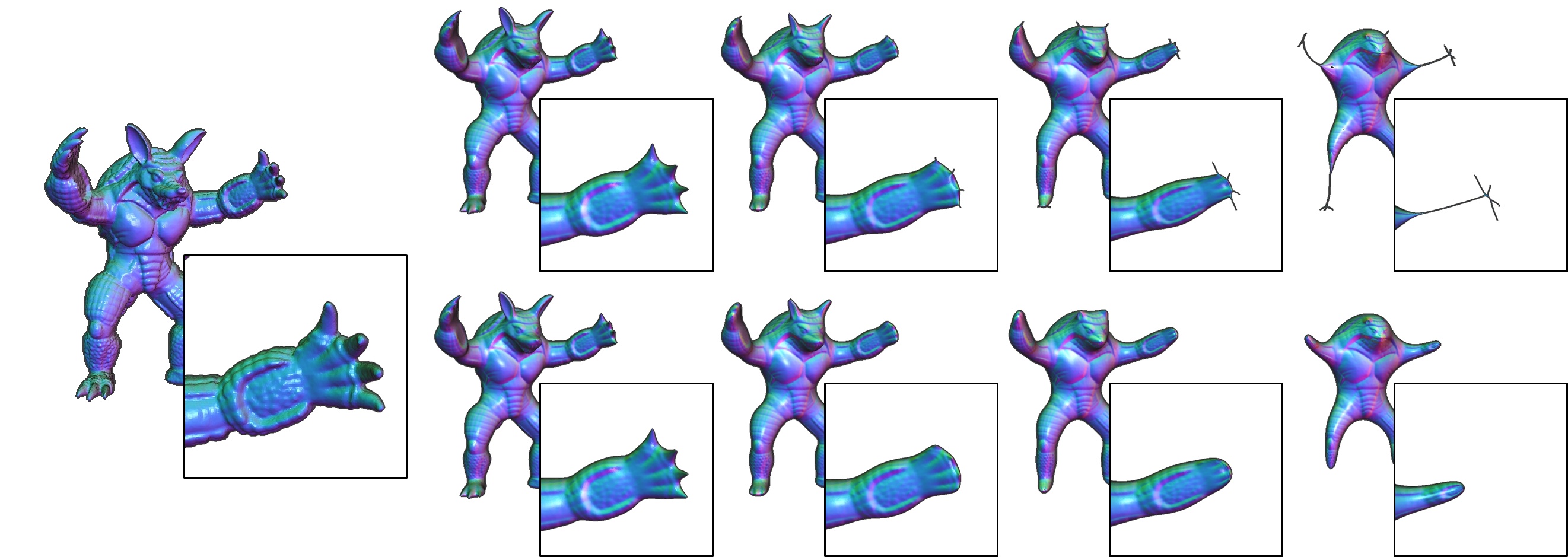}
  \caption{
  \label{f:teaser}
The armadillo man model (left) and the results of traditional MCF (top) compared to the results of our modified MCF. Both sets of results show the computed surface after 2, 5, 10, and 25 semi-implicit time-steps. The zoom-in on the left hand show that our modified flow successfully avoids forming neck pinching singularities that commonly occur in the mean-curvature flow of non-convex regions.
  }
}


\maketitle


\begin{abstract}
\input{abstract}

\end{abstract}

\begin{CRcatlist}
  \CRcat{I.3.5}{Computer Graphics}{Computational Geometry and Object Modeling}{Curve, surface, solid, and object representations}
\end{CRcatlist}

\keywordlist

\section{Introduction}
\label{s:intro}
\input{intro}

\section{Review of Mean-Curvature Flow}
\label{s:review}
\input{review}

\section{Modifying the Flow}
\label{s:modifying}
\input{modifying}

\section{Results and Discussion}
\label{s:results}
\input{results}

\section{Conclusion}
\label{s:conclusion}
\input{conclusion}

\ignore
{
\section{Notes}
\begin{itemize}
\item Yamabe and Ricci flows
\item Modify the introduction to highlight the fact that convergence seems to be conformal to a sphere
\item Foldovers?
\item Convergence? Specifically, after the shape of the surface has ``converged'' to a sphere, could there still be tangential drift?
\item Show plot of (some) convexity measure as a function of the number of time-steps.
\item Make the distinction between poor conditioning and singularities
\end{itemize}
}



\bibliographystyle{acmsiggraph}
\bibliography{paper}

\appendix
\section{Continuous Formulation}
\label{a:continuous}
\input{appendix_continuous}

\section{Analytic Flow Solutions}
\label{a:analytic}
\input{appendix_analysis}

\section{Energy of the Flow}
\label{a:energy}
\input{appendix_energy}

\end{document}

%% file: abstract.tex
This work considers the question of whether mean-curvature flow can be modified to avoid the formation of singularities. We analyze the finite-elements discretization and demonstrate why the original flow can result in numerical instability due to division by zero. We propose a variation on the flow that removes the numerical instability in the discretization and show that this modification results in a simpler expression for both the discretized and continuous formulations. We discuss the properties of the modified flow and present empirical evidence that not only does it define a stable surface evolution for genus-zero surfaces, but that the evolution converges to a conformal parameterization of the surface onto the sphere.

%% file: intro.tex
Mean-Curvature flow is one of the more basic flows that has been used to evolve surface geometry. It can be equivalently formulated as a flow that either (1)~minimizes the gradient of the surface embedding or (2)~minimizes surface area. As the former, it has played an essential role in the area of mesh fairing, removing noise by smoothing the embedding~\cite{Taubin:SIGGRAPH:1995,Desbrun:SIGGRAPH:1999}. As the latter, it has been essential in the study of minimal surfaces~\cite{Chopp:JCP:1993,Pinkall:EM:1993}.

Despite its pervasiveness, the utility of mean-curvature flow has been restricted by the formation of singularities during the course of the flow. As a result, convergence proofs have been limited to a class of simple shapes (e.g. mean-convex surfaces)~\cite{Huisken:JDG:1984}.

In this work, we propose a modification of traditional mean-curvature flow that may address this limitation.
We proceed by analyzing the finite-elements discretization of the flow, identify a potential for division-by-zero in the definition of the system matrix, and show that a minor modification to the flow removes the numerical instability,
providing simpler expressions for both the discrete and continuous formulations of the flow. An analysis of this flow shows that the modification does not change the flow in spherical regions, and it slows down the evolution in cylindrical regions, avoiding the formation of undesirable neck pinches. Though we do not have a proof that it will always happen, we show numerous examples of the flow on highly concave, genus-zero, surfaces, demonstrating that not only is the flow non-singular, but it also converges to a conformal parameterization of the surface onto a sphere.

An example of the modified mean-curvature flow can be seen in Figure~\ref{f:teaser}. The original Armadillo-Man model is shown on the left. The top row shows the results of traditional mean-curvature flow (with vertices of collapsed regions fixed and taken out of the linear system, following Au~{\em et al.}~\shortcite{Au:SIGGRAPH:2008}) while the bottom row shows results of the modified flow. Although both approaches smooth the geometry, the concavity of the shape results in singularities when evolving with the traditional flow. For example, the (cylindrical) fingers collapse within 5 iterations, forming neck pinches before they can be merged into the hand. In contrast, the modified flow avoids these types of singularities, slowing down the inward flow of the extremities and allowing them to merge into the appendages before they have an opportunity to collapse. The limit behavior of the flow can be seen in Figure~\ref{f:armadillo-conformality}, demonstrating that the flow evolves to a conformal parameterization of the model onto a sphere.


\subsection*{Outline}
The remainder of this paper is structured as follows. We present a review of both the continuous formulation of mean-curvature flow, and its discretization using finite-elements in Section~\ref{s:review}. We analyze the numerical stability of the discretization in Section~\ref{s:modifying} and show how the flow can be modified to avoid potential division-by-zero. We evaluate our approach in Section~\ref{s:results}, where we show the flow for a variety of genus-zero surfaces and discuss its properties. We conclude in Section~\ref{s:conclusion}, summarizing our work.

%% file: review.tex
In this section, we briefly review mean-curvature flow. We start with the continuous formulation and then describe a semi-implicit, finite-elements discretization. Although both are classic derivations (see e.g. \cite{mantegazza} for a modern treatment of the continuous formulation, and \cite{dziuk1990} for the finite-elements discretization), we repeat them here for completeness, as they are required for understanding the motivation for our modified flow. 

\subsection{Continuous Formulation}
Informally, mean-curvature flow can be thought of as a flow that pushes a point on a surface towards the average position of its neighbors. As with image filtering (when replacing a pixel's value by the average value of its neighbors) this flow has the effect of smoothing out the geometry.

\ifjakeversion
\noindent\textbf{Definition (MCF)}: Let $M$ be a two dimensional manifold, let $\Phi_t:M\rightarrow\R^3$ be a smooth family of immersions, and let $g_t(\cdot,\cdot)$ be the metric induced by the immersion at time $t$. We say that $\Phi_t$ is a solution to the {\em mean-curvature flow} if:
\begin{equation}
\frac{\partial\Phi_t}{\partial t} = \Delta_t\Phi_t
\label{eq:mcf}
\end{equation}
where $\Delta_t$ is the Laplace-Beltrami operator defined with respect to the metric $g_t$.
\else
Formally, given a manifold $M$ and given an embedding $\Phi:M\rightarrow\R^3$, mean-curvature flow can be characterized by a family of functions $\Phi_t:M\rightarrow{\mathbb R}^3$ that satisfy:
\begin{equation}
\frac{\partial\Phi_t}{\partial t}=\Delta_t \Phi_t\qquad\hbox{and}\qquad
\Phi_\tinyZ(p)= \Phi(p)\quad\forall p\in M
\label{eq:mcf}
\end{equation}
\ifjakeversion
where $\Delta_t$ is the Laplace-Beltrami operator defined with respect to the metric $g_t(\cdot,\cdot)$, induced by the embedding of $\Phi_t$.
\else
where $\Delta_t$ is the Laplace-Beltrami operator defined by $\Phi_t$.
\fi
\fi


\subsection{Finite Elements Discretization}
To model the flow in practice, we need to discretize the differential equation. We transform the continuous system of equations into a finite-dimensional system by choosing a function basis $\{B_1,\ldots,B_N\}:M\rightarrow\R$. Using this basis, we represent the map at time $t$ by the coefficient vector $\vec{x}(t)=\{x_1(t),\ldots,x_N(t)\}\subset\R^3$ so that:
$$\Phi_t(p) = \sum_{i=1}^N x_i(t) B_i(p).$$
Although we cannot solve Equation~\ref{eq:mcf} exactly, since the solution is not guaranteed to be within the span of the $\{B_i\}$, we can solve the system in a least-squares sense using the Galerkin formulation:
$$
\mathop{\int}_M\left(\frac{\partial\Phi_t}{\partial t} \cdot B_i\right)\measure{t} = \mathop{\int}_M\left( \Delta_t\Phi_t \cdot B_i\right)\measure{t}\quad\forall\ 1\leq i \leq N
$$
\ifjakeversion
with $\measure{t}$ the volume form induced by the metric $g_t$.
\else
with $\measure{t}$ the volume form induced by the embedding $\Phi_t$. These are in fact three separate equations, for each of the components of $\Phi_t$.
\fi

Setting $D^t$ and $L^t$ to be the mass and stiffness matrices of the embedding at time $t$:\footnote{We implicitly assume that either $M$ is water-tight or that its boundary is fixed throughout the course of the flow.}
$$
D^t_{ij}=  \mathop{\int}_M (B_i\cdot B_j)\measure{t}\qquad
L^t_{ij}=- \mathop{\int}_M g_t\left(\nabla_t B_i , \nabla_t B_j\right)\measure{t}
$$
a semi-implicit time discretization defines the linear system relating coefficients at time $t+\delta$ to coefficients at time $t$:
$$
\left(D^t-\delta L^t\right)\vec{x}(t+\delta)=D^t\vec{x}(t).
$$
\ifjakeversion
In the above equations, $\nabla_t$ denotes the gradient operator with respect to the metric $g_t$. 
\else
In the above equations, $\nabla_t$ denotes the gradient operator with respect to the embedding $\Phi_t(M)$ and $g_t(\cdot,\cdot)$ is the metric induced by the embedding $\Phi_t$.
\fi

%% file: modifying.tex
We begin this section by considering how numerical instabilities can arise with traditional mean-curvature flow and then propose a modified flow that resolves this problem.

\subsection{Numerical Instability}
\ifjakeversion
A challenge of using mean-curvature flow becomes apparent when we compute the coefficients of the matrix $D^t-\delta L^t$ by integrating with respect to the metric defined by the original embedding, $\Phi_\tinyZ$, rather than the current embedding, $\Phi_t$.
\else
A challenge of using mean-curvature flow becomes apparent when we compute the coefficients of the matrix, $D^t-\delta L^t$, by integrating with respect to the metric on the original surface, $\Phi_\tinyZ(M)$, rather than $\Phi_t(M)$.
\fi

\ifjakeversion
To this end, we consider how the geometry is stretched over the course of the flow, characterized by the endomorphism $g_\tinyZ^{-1}\cdot g_t$.\footnote{Here we identify $g_t$ with the map from the tangent space to its dual so that the map $g_\tinyZ^{-1}\cdot g_t$ is a well-defined map from the tangent space to itself.} This operator is self-adjoint with respect to both $g_\tinyZ$ and $g_t$. Its eigenvectors, $v_1$ and $v_2$, define the principal directions of stretch (orthogonal with respect to both $g_\tinyZ$ and $g_t$) and its eigenvalues, $\lambda_1^2$ and $\lambda_2^2$, define the (squares of the) magnitudes of stretch along these directions. 
\else
To this end, we consider how the map $\Phi_t\circ\Phi_\tinyZ^{-1}$ stretches the geometry, characterized by the endomorphism $g_\tinyZ^{-1} \cdot g_t$. This operator is self-adjoint with respect to both $g_\tinyZ$ and $g_t$. Its eigenvectors, $v_1$ and $v_2$, define the principal directions of stretch (orthogonal with respect to both $g_\tinyZ$ and $g_t$) and its eigenvalues, $\lambda_1^2$ and $\lambda_2^2$, define the (squares of the) magnitudes of stretch along these directions. 
\fi

Note that both the stretch directions, $v_i$, and the stretch factors, $\lambda_i$, depend on the time parameter $t$. However, we omit it in our notation for simplicity.




\subsubsection*{The Mass Matrix}
Using the chain rule, we obtain an expression for the $(i,j)$-th coefficients of the mass matrix as:
\ifcleanratio
$$D_{ij}^t = \mathop{\int}_M B_i \cdot B_j\cdot\sqrt{\left|\dg\right|}\measure{\tinyZ},$$
where $\sqrt{\left|\dg\right|}$ gives the ratio of area elements. Since mean-curvature flow is area minimizing, the values of $\left|\dg\right|$ tend to be small so the computation of $D_{ij}^t$ is numerically stable.
\fi

\subsubsection*{The Stiffness Matrix}
The situation gets more complicated when computing the stiffness matrix. The challenge here is due to the fact that, as area shrinks, the corresponding derivatives grow. As a result, since mean-curvature flow is area minimizing, there is potential for the Laplacian to blow up.

To make this explicit, we decompose the gradient in the integrand of $L_{ij}^t$ into orthogonal components along the principal directions of stretch, $v_1$ and $v_2$. Then, using the fact that shrinking the domain of a function by a factor of $\lambda$ scales its derivative by the same factor, the expression for the $(i,j)$-th coefficient of the stiffness matrix becomes:
\ifcleanratio
\begin{align}
\nonumber
L_{ij}^t
&= -\mathop{\int}_M\sum_{k,l=1}^2\frac{\partial B_i}{\partial v_k}\frac{\partial B_j}{\partial v_l}g_t^{k,l}\sqrt{\left|\dg\right|}\measure{\tinyZ}\\
\nonumber
&= -\mathop{\int}_M\sum_{k=1}^2\frac{\partial B_i}{\partial v_k}\frac{\partial B_j}{\partial v_k}\frac{1}{g_t(v_k,v_k)}\lambda_1\lambda_2\measure{\tinyZ}\\
\label{eq:lap-instability}
&= -\mathop{\int}_M\sum_{k=1}^2\frac{\partial B_i}{\partial v_k}\frac{\partial B_j}{\partial v_k}\frac{1}{g_\tinyZ(v_k,v_k)}\frac{\lambda_1\lambda_2}{\lambda_k\lambda_k}\measure{\tinyZ}
\end{align}
\fi
where $g_t^{k,l}$ is the $(k,l)$-th coefficient of the inverse of the $2\times2$ matrix whose $(i,j)$-th entry is $g_t(v_i,v_j)$.

Examining Equation~\ref{eq:lap-instability}, we observe that while the value $\frac{\partial B_i}{\partial v_k}\frac{\partial B_j}{\partial v_k}\frac{1}{g_\tinyZ(v_k,v_k)}$ is stable, as it only depends on the values of the partial derivatives of the $B_i$ along directions that are unit-length under the initial metric, $g_\tinyZ$, the stretch ratios $\lambda_1\lambda_2/\lambda_2\lambda_2=\lambda_1/\lambda_2$ and $\lambda_1\lambda_2/\lambda_1\lambda_1=\lambda_2/\lambda_1$ might not be.
\ifjakeversion
In particular, as the stretching becomes more anisotropic (i.e. the mapping $\Phi_t$ becomes less and less conformal with respect to the metric $g_\tinyZ$) one of the two ratios will tend to infinity and the integrand will blow up.
\else
In particular, as the stretching of $\Phi_\tinyZ(M)$ onto $\Phi_t(M)$ becomes more anisotropic (i.e. the mapping $\Phi_t\circ\Phi_\tinyZ^{-1}$ becomes less and less conformal) one of the two ratios will tend to infinity and the integrand will blow up.
\fi

\subsection{Conformalizing the Metric}
Since the instability of the stiffness matrix is the result of anisotropic stretching, we address this problem by replacing the metric $g_t$ in the computation of the system coefficients by the closest metric that is conformal to $g_\tinyZ$:
\ifcleanratio
$$\tilde{g}_t = \sqrt{\left|\dg\right|}g_\tinyZ.$$
\fi

\subsubsection*{FEM Discretization}
For the finite-elements discretization, the implementation of this modification is trivial. Instead of requiring that the stiffness matrix be computed anew at each time-step (as is required by traditional mean-curvature flow), the modified flow simply re-uses the stiffness matrix from time $t=0$.

Specifically, since the conformalized metric $\tilde{g}_t$ has the same determinant as the old metric $g_t$, the coefficients of the mass matrix are the same regardless of which of the two metrics we use. However, because the two eigenvalues of the conformalized metric are equal, $\tilde{\lambda}_1=\tilde{\lambda}_2$, the coefficients of the modified stiffness matrix become:
\begin{align*}
\tilde{L}_{ij}^t
&= -\int\limits_M\sum_{k=1}^2\frac{\partial B_i}{\partial v_k}\frac{\partial B_j}{\partial v_k}\frac{1}{g_\tinyZ(v_k,v_k)}\frac{\tilde{\lambda}_1\tilde{\lambda}_2}{\tilde{\lambda}_k\tilde{\lambda}_k}\measure{\tinyZ}\\
&= -\int\limits_M \sum_{k,l=1}^2\frac{\partial B_i}{\partial v_k}\frac{\partial B_j}{\partial v_l}g_\tinyZ^{k,l}\measure{\tinyZ} =
L_{ij}^0,
\end{align*}
making it independent of time $t$.

\subsubsection*{Continuous Formulation}
\ifjakeversion
The conformalization of the metric also results in a simpler continuous formulation, allowing us to replace the PDE in Equation~\ref{eq:mcf} with the following.

\noindent\textbf{Definition (cMCF)}: Let $(M,h)$ be a two dimensional Riemannian manifold (with metric $h$), let $\Phi_t:M\rightarrow\R^3$ be a smooth family of immersions, and let $g_t(\cdot,\cdot)$ be the metric induced by the immersion at time $t$. We say that $\Phi_t$ is a solution to the {\em conformalized mean-curvature flow} if:
\begin{equation}
\frac{\partial\Phi_t}{\partial t} = \sqrt{\left|\dgh\right|}\Delta_h\Phi_t
\label{eq:cmcf}
\end{equation}
where $\Delta_h$ is the Laplace-Beltrami operator defined with respect to the metric $h$.
\else
The conformalization of the metric also results in a simpler formulation of the continuous flow, allowing us to replace the PDE defining traditional mean-curvature flow:
\begin{equation}
\frac{\partial\Phi_t}{\partial t} = \Delta_t\Phi_t
\label{eq:continuous}
\end{equation}
which requires computing the new Laplace-Beltrami operator at each time $t$, with the simpler PDE:
\ifcleanratio
\begin{equation}
\frac{\partial\Phi_t}{\partial t} = \sqrt{\left|\dgi\right|}\Delta_\tinyZ\Phi_t.
\label{eq:continuous-modification}
\end{equation}
\else
\begin{equation}
\frac{\partial\Phi_t}{\partial t} = \sqrt{\frac{|g_\tinyZ|}{|g_t|}}\Delta_\tinyZ\Phi_t.
\label{eq:continuous-modification}
\end{equation}
\fi
\fi
(For a derivation, see Appendix~\ref{a:continuous}.)

\ifjakeversion
Note that if the mapping $\Phi_t$ is conformal with respect to $h$, the Laplace-Beltrami operators $\Delta_h$ and $\Delta_t$ are related by:
$$\Delta_t = \sqrt{\left|\dgh\right|}\Delta_h$$
and the flows in Equations~\ref{eq:mcf} and~\ref{eq:cmcf} are the same.
\else
Note that when the mapping $\Phi_t\circ\Phi_\tinyZ^{-1}$ is conformal, the Laplace-Beltrami operators $\Delta_\tinyZ$ and $\Delta_t$ are related by:
\ifcleanratio
$$\Delta_t = \sqrt{\left|\dgi\right|}\Delta_\tinyZ$$
\else
$$\sqrt{|g_t|}\Delta_t = \sqrt{|g_\tinyZ|}\Delta_\tinyZ$$
\fi
and the flows in Equations~\ref{eq:continuous} and~\ref{eq:continuous-modification} are the same.
\fi

%% file: results.tex
\ifjakeversion
We begin by examining some examples of the conformalized mean-curvature flow and then proceed to a discussion of its properties. In our discussion, we assume that an initial embedding $\Phi_\tinyZ:M\rightarrow\R^3$ is given, and we take $h$ to be the metric induced by this embedding, $h=g_\tinyZ$.
\else
We begin by examining some examples of the modified mean-curvature flow and then proceed to a discussion of its properties.
\fi

\subsection{Flowing Surfaces}
\ifjakeversion
To better understand how the conformalized flow evolves the embedding, we compare with two other flows. The first is traditional mean-curvature flow, which updates both the mass- and stiffness-matrix at each time-step. The second is the simple heat flow with respect to the metric $h=g_\tinyZ$ that keeps both matrices fixed, and has been proposed for efficient short-term flows in the context of mesh-fairing applications~\cite{Desbrun:SIGGRAPH:1999}. 
\else
To better understand the properties of the modified flow, we compare with two other flows. The first is traditional mean-curvature flow, which updates both the mass- and stiffness-matrix at each time-step. The second is the simple heat flow with respect to the metric $g_\tinyZ$ that keeps both matrices fixed, and has been proposed for efficient short-term flows in the context of mesh-fairing applications~\cite{Desbrun:SIGGRAPH:1999}. 
\fi

\subsubsection*{Analytic Flow}
We start by considering three simple examples for which an analytic expression of the flow can be computed. Understanding these simple cases provides some intuition as to why the modified flow might be free of singularities. These examples include the (hyperbolic) catenoid, the (elliptic) sphere, and the (parabolic) infinite cylinder. For each of these geometries and all three flows, the evolved surfaces can be characterized by their radius, $r(t)$, as described in Table~\ref{t:radii}. (For a derivation, see Appendix~\ref{a:analytic}.)

\ifjakeversion
\begin{table}[!t]
\centering
\setlength{\tabcolsep}{2.4pt}%
\renewcommand{\arraystretch}{1.2}%
{\footnotesize%
\begin{tabular}{l|l|l|l|}
 & \multicolumn{1}{c|}{MCF} & \multicolumn{1}{c|}{Heat Flow} & \multicolumn{1}{c|}{cMCF} \\
\hline
Catenoid & $r(t)=1$           & $r(t)=1$       & $r(t)=1$           \\
\hline
Sphere   & $r(t)=\sqrt{1-4t}$ & $r(t)=e^{-2t}$ & $r(t)=\sqrt{1-4t}$ \\
\hline
Cylinder & $r(t)=\sqrt{1-2t}$ & $r(t)=e^{-t}$  & $r(t)=1-t$         \\
\hline
\end{tabular}
}%
\caption{
Radii of the catenoid, sphere, and infinite cylinder under the different flows, as a function of time $t$.
\label{t:radii}
}
\vspace{-.1in}
\end{table}
\else
\begin{table}[!t]
\centering
\setlength{\tabcolsep}{2.4pt}%
\renewcommand{\arraystretch}{1.2}%
{\footnotesize%
\begin{tabular}{l|l|l|l|}
 & \multicolumn{1}{c|}{Traditional MCF} & \multicolumn{1}{c|}{Heat Flow} & \multicolumn{1}{c|}{Modified MCF} \\
\hline
Catenoid & $\,\,\,\quad r(t)=1$           & $r(t)=1$       & $\quad r(t)=1$           \\
\hline
Sphere   & $\,\,\,\quad r(t)=\sqrt{1-4t}$ & $r(t)=e^{-2t}$ & $\quad r(t)=\sqrt{1-4t}$ \\
\hline
Cylinder & $\,\,\,\quad r(t)=\sqrt{1-2t}$ & $r(t)=e^{-t}$  & $\quad r(t)=1-t$         \\
\hline
\end{tabular}
}%
\caption{
Radii of the catenoid, sphere, and infinite cylinder under the different flows, as a function of time $t$.
\label{t:radii}
}
\vspace{-.1in}
\end{table}
\fi

Examining this table, we make several observations.

First, for each shape, the derivatives of the three flows at $t=0$ are equal, since all three start by flowing the embedding along the normal direction, with speed equal to the negative of the mean-curvature.

Second, the catenoid remains fixed under all three flows since its mean-curvature is everywhere zero.

Third, the heat flow remains stable for all shapes. This is expected since computing the flow is equivalent to repeatedly multiplying by the inverse of the system matrix (giving the characteristic exponent in the radius function) so the long term-behavior can be computed by projecting onto the lower eigenvectors of the system. 

\ifjakeversion
Fourth, for the case of the sphere, both MCF and cMCF give the same result. This is because the mean-curvature flow of the embedding of the sphere is conformal, so we have $\tilde{g}_t=g_t$, and both flows define the same linear system.
\else
Fourth, for the case of the sphere, both the traditional and the modified mean-curvature flow give the same result. This is because the mean-curvature flow of the sphere is conformal, so we have $\tilde{g}_t=g_t$, and both the flows define the same linear system.
\fi

\begin{figure}[h]
\center
{
\includegraphics[width=0.95\columnwidth]{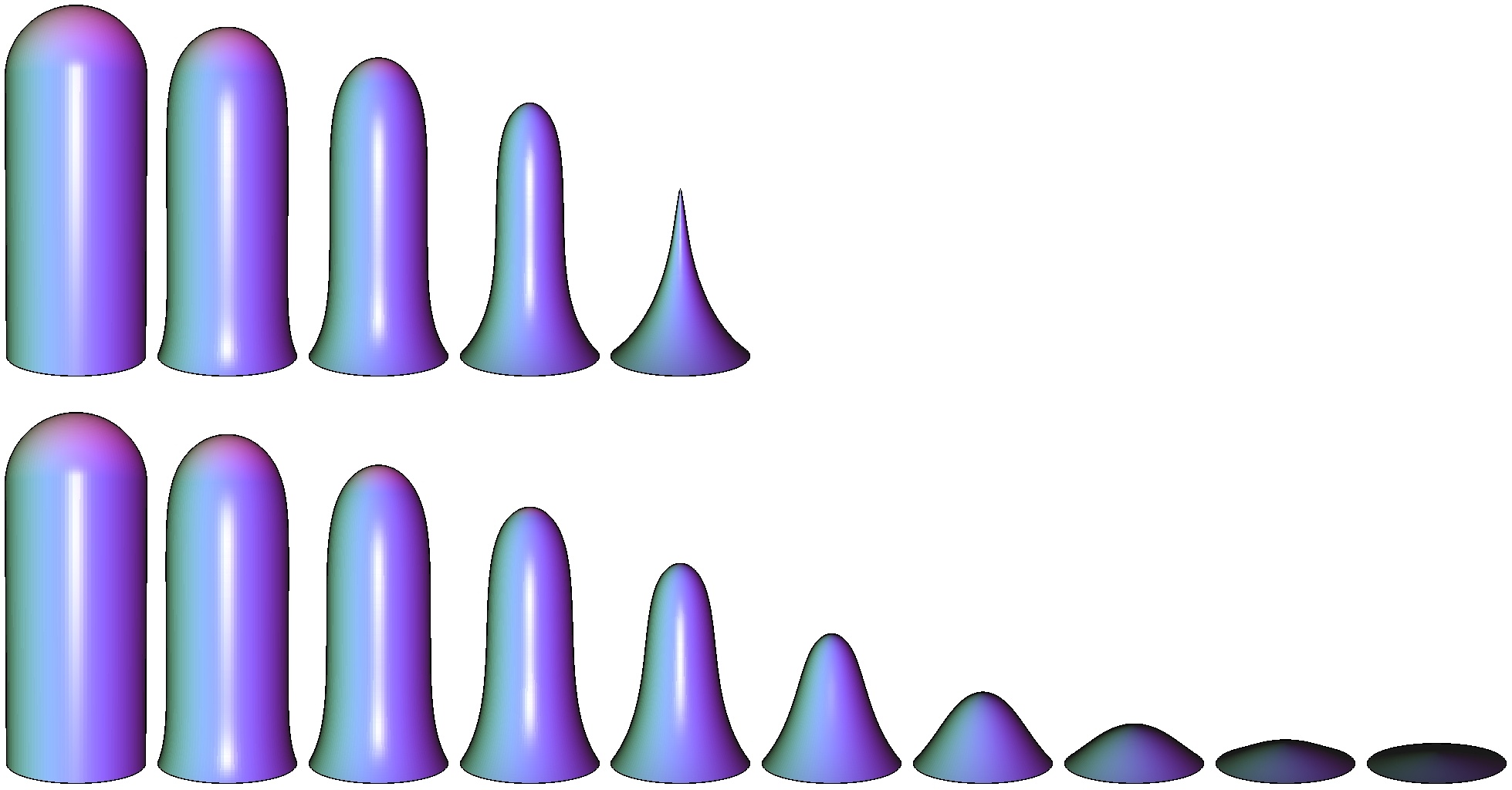}
}
  \caption{
  \label{f:short_cyl_p}
Evolution of an embedding of a rounded cylinder using MCF (top) and the modified flow (bottom). Note that with cMCF the cylinder collapses more slowly, allowing the spherical cap to ``catch-up'', avoiding the singularity.
  }
\end{figure}

\ifjakeversion
Finally, for the cylinder, cMCF slows the rate of shrinking so that the flow towards the cylinder's axis no longer accelerates with time. This is demonstrated in Figure~\ref{f:short_cyl_p}, where the first few iterations of MCF and cMCF are shown. We can see that flowing with cMCF, the cylindrical center collapses more slowly, allowing the spherical top to ``catch-up'' avoiding the formation of the singularity that appears in MCF. 
\else
Finally, for the cylinder, the modified mean-curvature flow slows the rate of shrinking so that the flow towards the cylinder's axis no longer accelerates with time. This is demonstrated in Figure~\ref{f:short_cyl_p}, where the first few iterations of MCF and the modified flow are shown. We can see that the cylindrical center collapses more slowly, allowing the spherical top to ``catch-up'' avoiding the singularity that occurs in MCF. 
\fi

As we will see next, the same effect allows surface extremities to collapse into the main body without forming neck pinches and, for genus-zero surfaces, evolves into a conformal parameterization of the surface onto a sphere.

\begin{figure*}[!t]
\center
{
	\includegraphics[width=2\columnwidth]{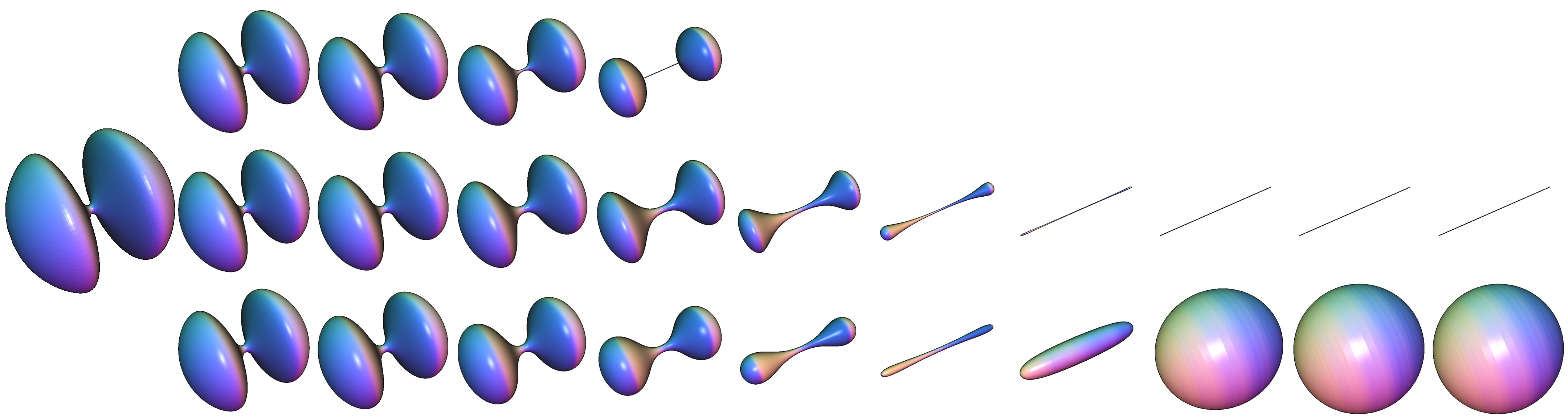}
}
  \caption{
  \label{f:dumbbell}
Evolution of a dumbbell model with step-size $\delta=1\times10^{-3}$.
  }
\end{figure*}

\begin{figure*}[!t]
\center
{
	\includegraphics[width=2\columnwidth]{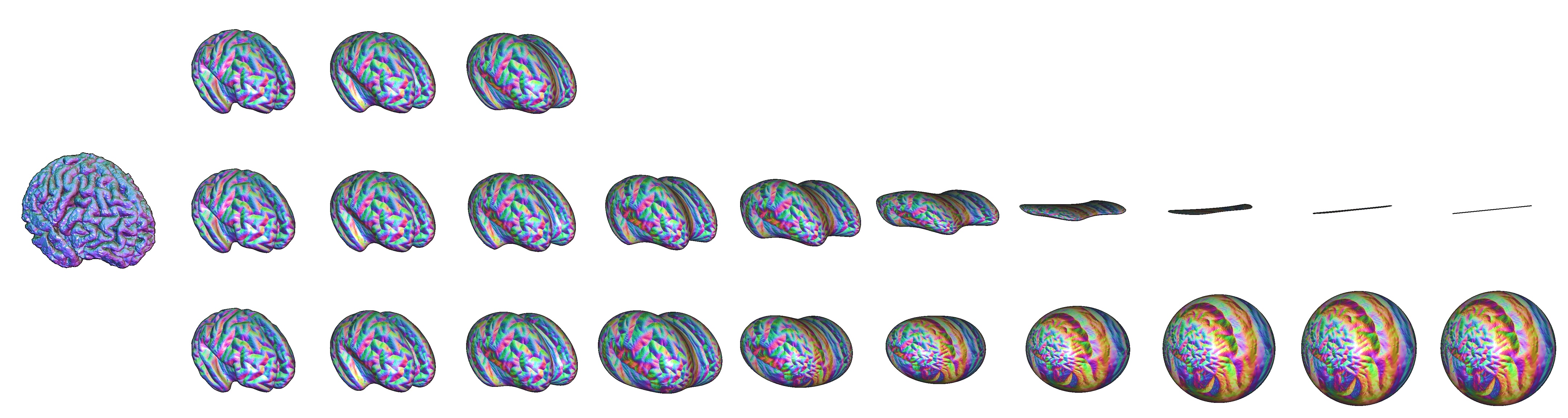}
}
  \caption{
  \label{f:brain}
Evolution of brain gray-matter with step-size $\delta=5\times10^{-4}$.
 }
\end{figure*}

\begin{figure*}[!t]
\center
{
	\includegraphics[width=2\columnwidth]{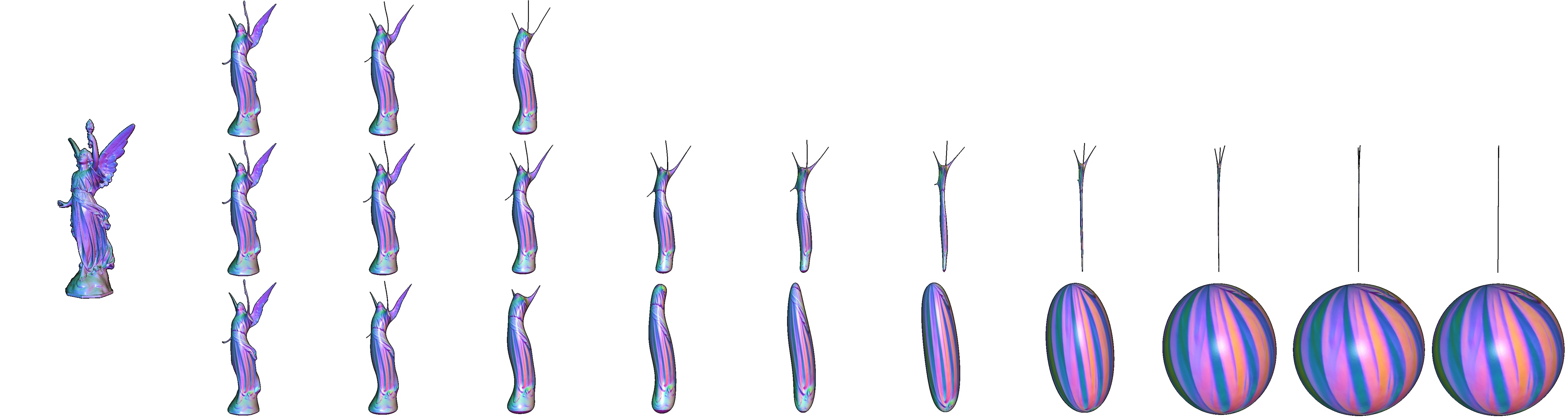}
}
  \caption{
  \label{f:lucy}
Evolution of the Lucy model with step-size $\delta=5\times10^{-4}$.
  }
\end{figure*}

\iffullimages
	\begin{figure*}[!t]
	\center
	{
		\includegraphics[width=2\columnwidth]{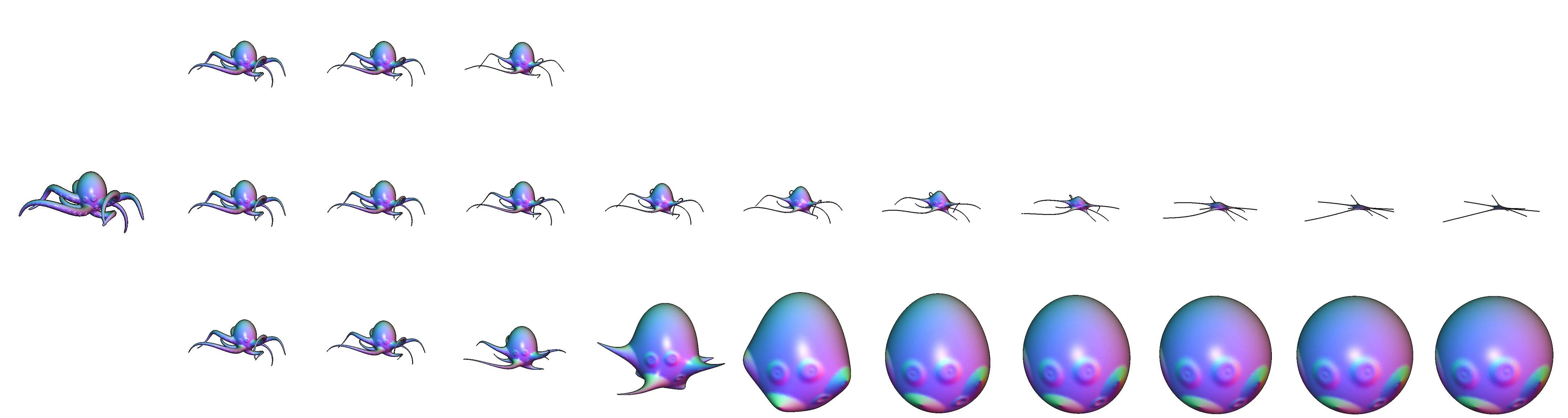}
	}
 	 \caption{
 	 \label{f:octopus}
	Evolution of an octopus model with step-size $\delta=2\times10^{-4}$.
	  }
	\end{figure*}
\fi

\begin{figure*}[!t]
\center
{
	\includegraphics[width=2\columnwidth]{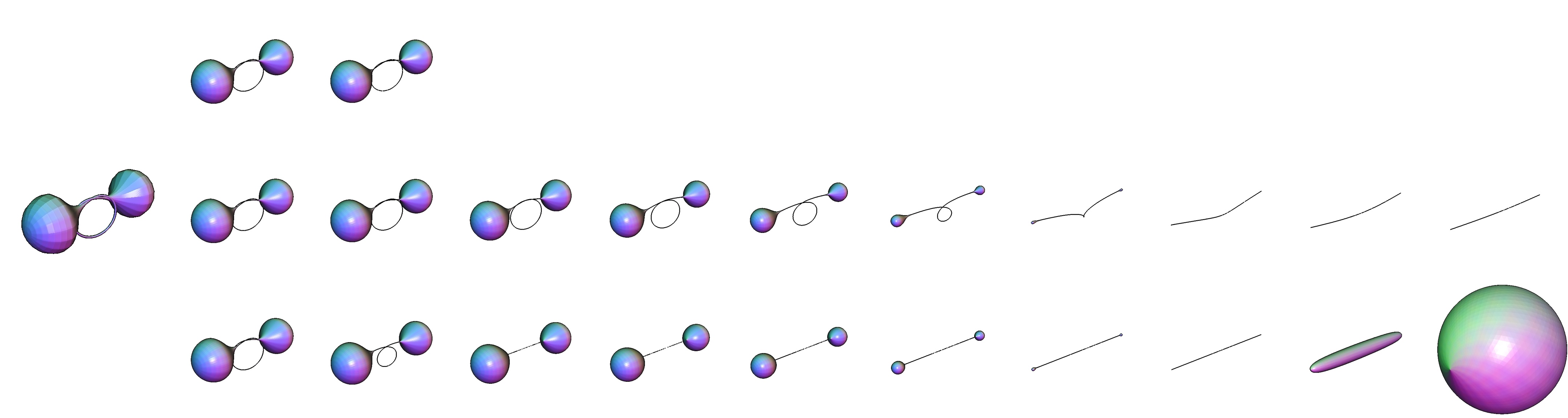}
}
  \caption{
  \label{f:knot}
Evolution of a knot with step-size $\delta=5\times10^{-4}$.
  }
\end{figure*}

\iffullimages
	\begin{figure*}[!t]
	\center
	{
		\includegraphics[width=2\columnwidth]{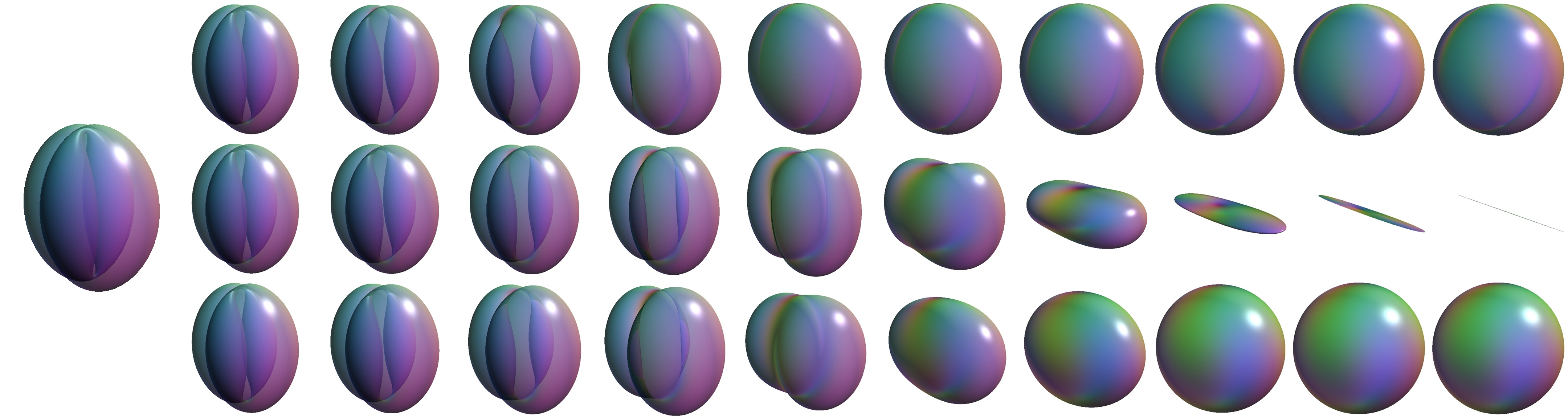}
	}
	  \caption{
	  \label{f:ds01}
	Evolution of the surface of an eigenfunction of the Dirac operator with step-size $\delta=5\times10^{-4}$.
	}
	\end{figure*}
\fi





\ignore
{
\begin{figure*}[!t]
\center
{
	\includegraphics[width=2\columnwidth]{cylinder.jpg}
}
  \caption{
  \label{f:cylinder}
Evolution of a cylinder, with boundaries, with step-size $\delta=1\times10^{-3}$.
 }
\end{figure*}

\begin{figure*}[!t]
\center
{
	\includegraphics[width=2\columnwidth]{hand.jpg}
}
  \caption{
  \label{f:hand}
Evolution of a hand, with boundary, with step-size $\delta=1\times10^{-4}$.
 }
\end{figure*}
}

\subsubsection*{Empirical Evaluation}
We ran the modified flow on a number of genus-zero models. For each one, we defined the mass- and stiffness-matrices using the hat basis~\cite{Dziuk:LNM:1988}:
\begin{eqnarray*}
D_{ij}&=&
\left\{\begin{array}{ll}
\displaystyle\frac{|T_{ij}^1|+|T_{ij}^2|}{12} & \hbox{if }j\in N(i) \\
\displaystyle\sum_{k\in N(i)} D_{ik} & \hbox{if } j=i
\end{array}\right.\\
L_{ij}&=&
\left\{\begin{array}{ll}
\displaystyle\frac{\cot\beta_{ij}^1+\cot\beta_{ij}^2}{2} & \hbox{if }j\in N(i) \\
\displaystyle-\sum_{k\in N(i)} L_{ik} & \hbox{if } j=i
\end{array}\right.
\end{eqnarray*}
where $N(i)$ are the indices of the vertices adjacent to vertex $i$, $T_{ij}^1$ and $T_{ij}^2$ are the two triangles sharing edge $(i,j)$, and $\beta_{ij}^1$ and $\beta_{ij}^2$ are the two angles opposite edge $(i,j)$.

We performed the semi-implicit time-stepping using a direct CHOLMOD solver~\cite{David:JMAA:1999}, running for 512 time-steps, terminating early if numerical instabilities were identified.\footnote{Numerical instability was defined by failure of CHOLMOD to produce a solution, due to the fact that linear system was not positive definite.} Following~\cite{Huisken:JDG:1984}, we uniformly scaled the map after each step to obtain a surface with unit area. (This is equivalent to reducing the time-step size at larger values of $t$, providing a finer-grained sampling of the flow when the surface evolves more quickly.)

Our visualizations show the results of traditional mean-curvature flow (top), the heat flow (middle), and the modified flow (bottom). They show the original model on the left and the results of the flow for $2^0,\ldots,2^9$ time-steps to the right, with per-vertex colors assigned using the normals from the original surface.

Figure~\ref{f:dumbbell} shows an example of the flow for a dumbbell shape. Because of the concavity at the center, traditional mean-curvature flow quickly creates a singularity (before the 16-th time-step with $\delta=10^{-3}$) and the flow cannot proceed. In contrast, the modified flow slows down the collapse near the center, allowing the extremities to collapse more quickly, evolving the embedding into a narrow ellipse which then flows back to a map onto the sphere. 

While the heat flow remains stable, it does not tend to evolve towards a smooth shape. (We have validated the long-term behavior of this flow more formally by projecting the embedding function onto the lower-frequency eigenvectors of the Laplace-Beltrami operator.)

\iffullimages
Figures~\ref{f:brain}-\ref{f:ds01} show challenging examples of flow for highly non-convex shapes including brain gray-matter, the Lucy model, an octopus, a twisted dumbbell, and the self-intersecting surface of an eigenfunction of the Dirac operator~\cite{Crane:SIGGRAPH:2011}. Again, we find that MCF quickly creates a singularity and the flow cannot proceed. Similarly, the heat flow does not evolve towards a smooth embedding. It is only using cMCF that the embeddings evolve to maps to a sphere.
\ifsupplemental
\ifjakeversion
(For additional examples of the evolution of embeddings of genus-zero models under cMCF, please see the supplemental material.)
\else
(For additional examples of the evolution of embeddings of genus-zero models under the modified mean-curvature flow, please see the supplemental material.)
\fi
\fi
\fi

\subsection{Discussion}
Several questions arise when considering the evolution given by the modified (area-normalized) mean-curvature flow.

\iffullimages
\begin{figure*}
	\center{\includegraphics[width=2\columnwidth]{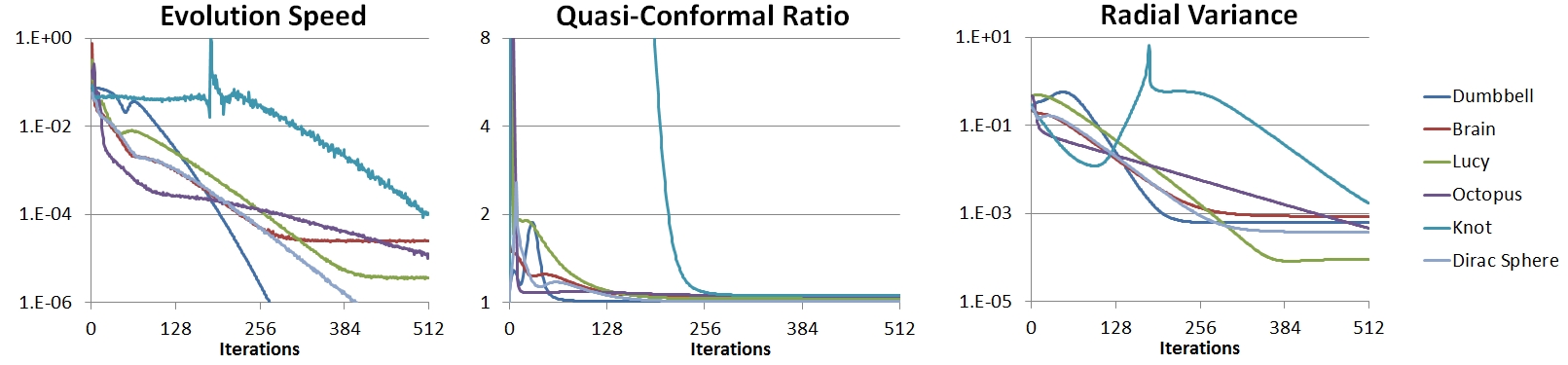}}
  \caption{
  \label{f:measures}
	Convergence (left), conformality (middle), and sphericity (right) of the modified mean-curvature flow for the models in Figures~\ref{f:dumbbell}-\ref{f:ds01}.
	}
\end{figure*}
\fi

\subsubsection*{Does it converge?}
\ifjakeversion
While this work is motivated by the goal of developing a variation of mean-curvature flow that is non-singular and converges when applied to embeddings of genus-zero surfaces, we have only been able to provide experimental confirmation of this property and leave the proof (or the existence of a counter-example) as an open question. Figure~\ref{f:measures} (left) empirically confirms the convergence of the flow for the models in Figure~\ref{f:dumbbell}-\ref{f:knot}, giving the magnitude of the difference between successive maps at each iteration.
\else
While this work is motivated by the goal of developing a variation of mean-curvature flow that is non-singular and converges when applied to genus-zero surfaces, we have only been able to provide experimental confirmation of this property and leave the proof (or the existence of a counter-example) as an open question.
Figure~\ref{f:measures} (left) empirically confirms the convergence of the flow for the models in Figure~\ref{f:dumbbell}-\ref{f:knot}, giving the magnitude of the change in the coordinate function at each iteration.
\fi

\ifconvergencediscussion
\subsubsection*{What does it converge to?}
\iffullimages
Looking at Figures~\ref{f:dumbbell}-\ref{f:ds01},
\else
Looking at Figures~\ref{f:dumbbell}-\ref{f:knot},
\fi
\ifjakeversion
we observe that modified mean-curvature flow appears to always converge to a map onto the sphere.
\else
we observe that modified mean-curvature flow appears to always converge to a sphere.
\fi
This is confirmed empirically in Figure~\ref{f:measures} (right), which plots the variance of the distance of the mesh vertices from its barycenter, as a function of the number of iterations. The plots shows that even though the flow may initially make the embedding less spherical, the variance decays in the limit.

\subsubsection*{How does it converge?}
We are also interested in characterizing the mapping from the original surface to the fixed point of the flow. Comparing the triangulation on the original surface to the triangulation on the limit surface (Figure~\ref{f:armadillo-conformality}), we see that the limit surface appears to preserve the aspect ratio of the triangles, suggesting that the mapping to the limit surface is conformal.

\ifsgpversion
We confirm this empirically by measuring the quasi-conformal error, computed as the area-weighted average of the ratios of the largest to smallest singular values of the map's Jacobian~\cite{Sander:SIGGRAPH:2001}. For the model in Figure~\ref{f:armadillo-conformality}, the average quasi-conformal error is $1.034$, which is 	comparable to the error of $1.033$ for the conformal spherical \begin{wrapfigure}{r}{45mm}
\vspace{-3mm}
\hspace{-8mm}
{} 
	\includegraphics[width=51mm]{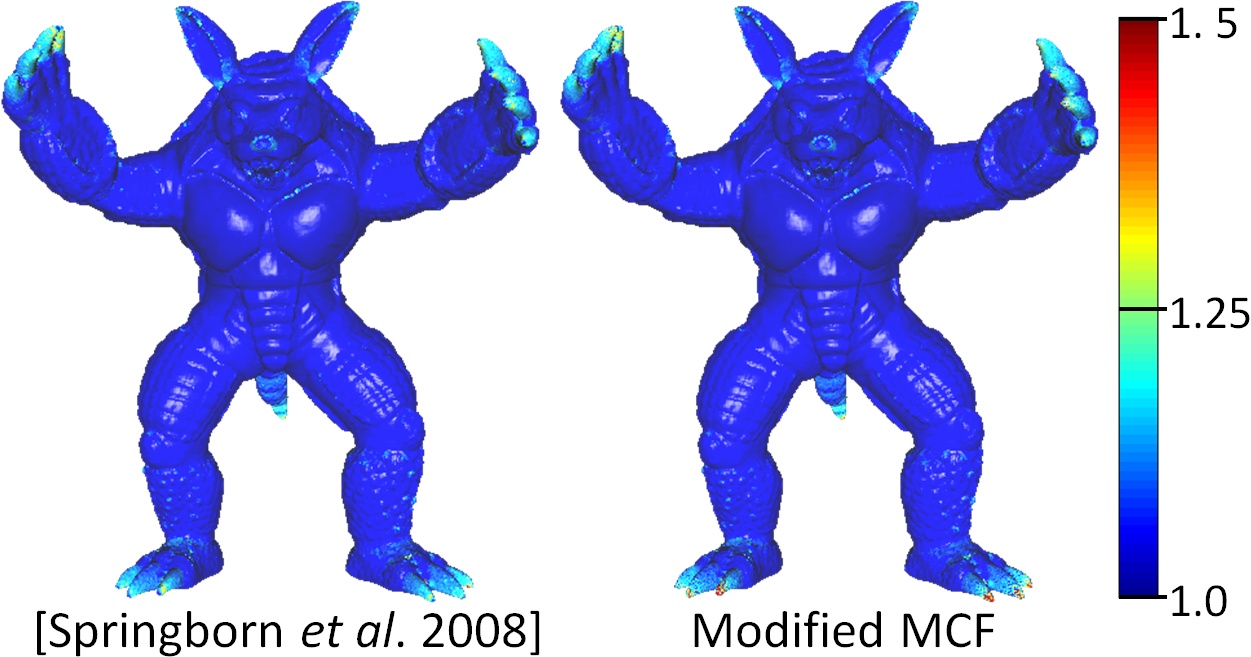}
\vspace{-6mm}
\end{wrapfigure}
parameterization of Springborn~{\em et al.}~\shortcite{Springborn:SIGGRAPH:2008}.
Visualizations of the errors for both of these
maps are shown in the inset on the right, highlighting the fact that the quasi-conformal errors are similarly distributed over the surface.
\else
We confirm this empirically by measuring the quasi-conformal error, computed as the area-weighted average of the ratios of the largest to smallest singular values of the mapping's Jacobian~\cite{Sander:SIGGRAPH:2001}. For the model in Figure~\ref{f:armadillo-conformality}, the average quasi-conformal error is $1.034$, which is comparable to the error of $1.033$ for the conformal spherical parameterization of
\begin{wrapfigure}{r}{45mm}
\vspace{-3mm}
\hspace{-8mm}
	\includegraphics[width=51mm]{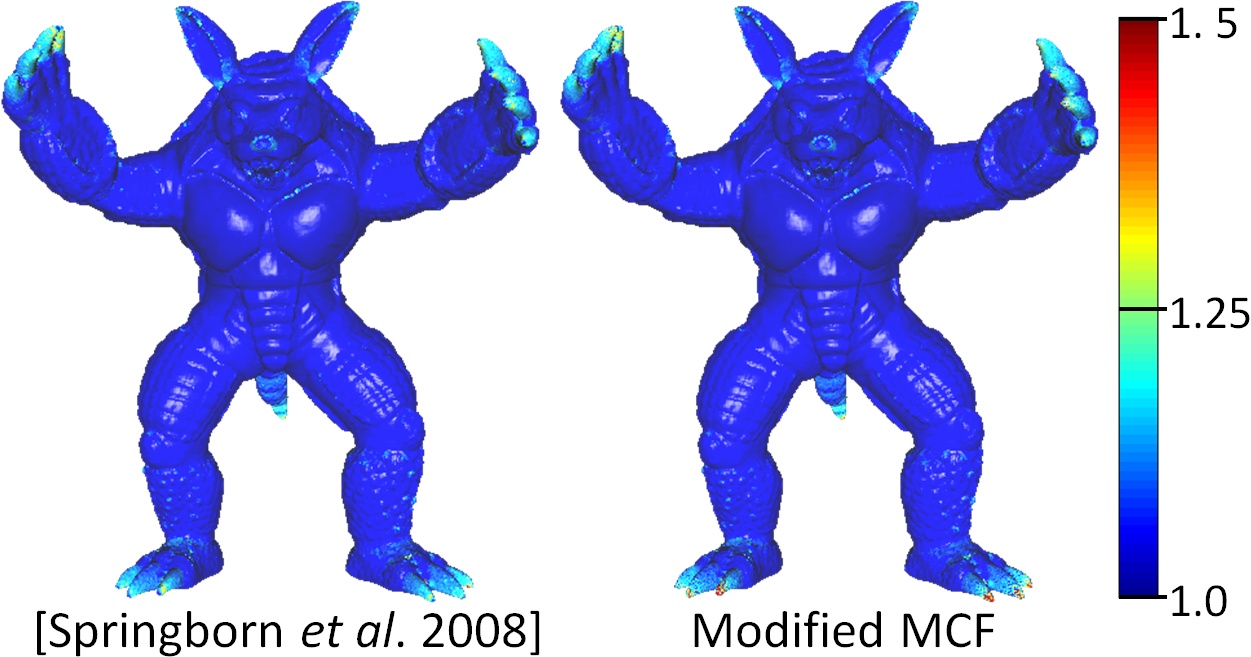}
\vspace{-6mm}
\end{wrapfigure}
Springborn~{\em et al.}~\shortcite{Springborn:SIGGRAPH:2008}. Visualizations of the errors for both of these maps are shown in the inset on the right, highlighting the fact that the quasi-conformal errors are similarly distributed over the surface.
\fi

\begin{figure}
	\center{\includegraphics[width=.8\columnwidth]{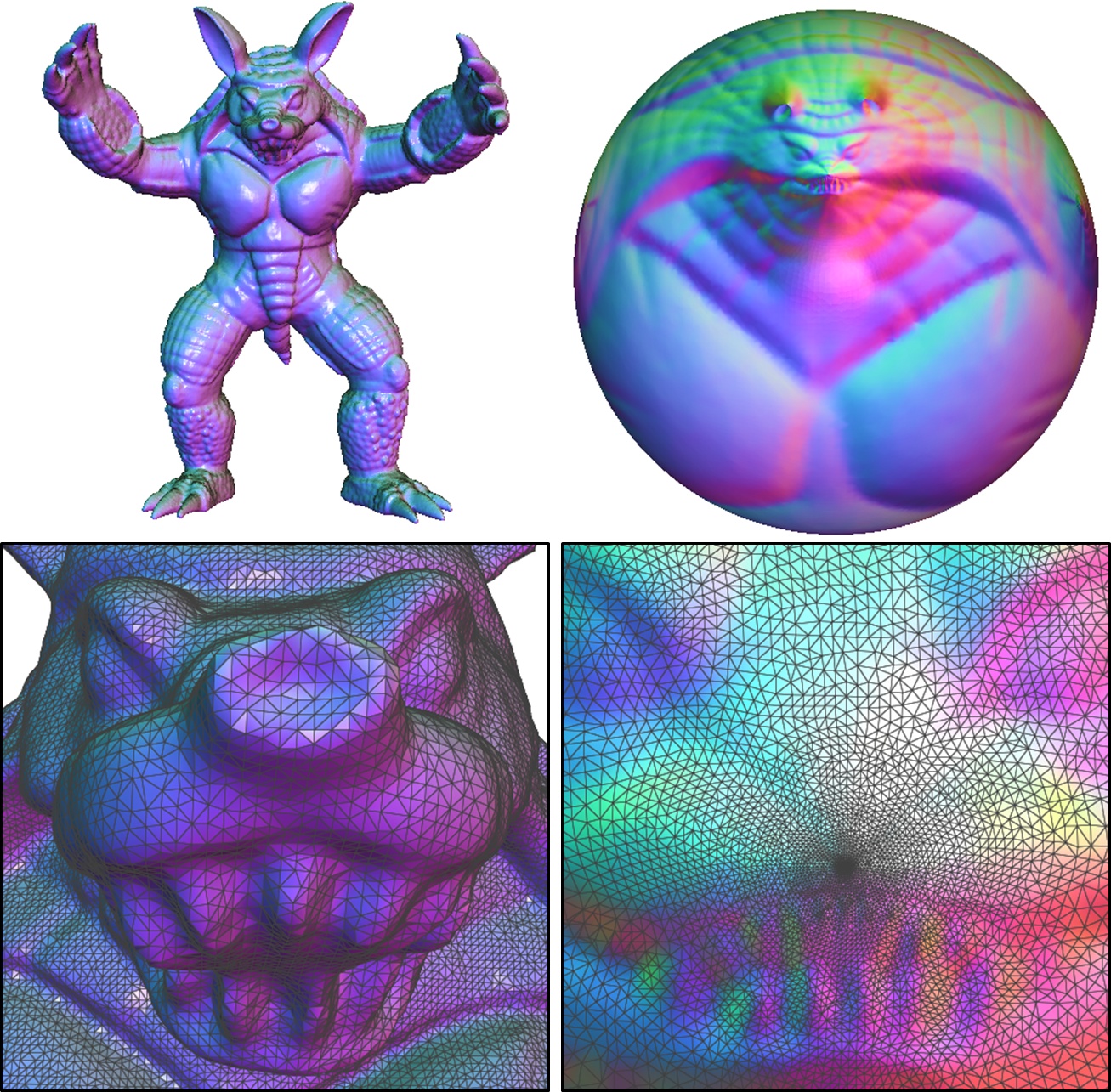}}
  \caption{
  \label{f:armadillo-conformality}
  The original armadillo-man model (left) and the surface obtained as the limit of the modified mean-curvature flow (right). Zoom-ins on the triangulation show that the mapping appears to preserve the aspect-ratio of the triangles, suggesting that the mapping is conformal.
	}
\end{figure}

More generally, Figure~\ref{f:measures} (middle) shows the plots of these ratios for the shapes in Figures~\ref{f:dumbbell}-\ref{f:knot} as a function of the number of iterations. We see that although the flow is not conformal, since the in-between maps have a high quasi-conformal error, the evolution does appear to converge to a conformal map, for all the shapes.

\noindent\textbf{Proposition}:
\newline\noindent
\ifjakeversion
If cMCF converges, than it converges to a map onto the sphere if and only if the limit map is conformal.
\else
If the modified mean-curvature flow converges, than it converges to a sphere if and only if the limit map is conformal.
\fi

\noindent\textbf{Proof}:
\newline\noindent
\ifjakeversion
($\Leftarrow$) If the map $\Phi_{t^*}$ is conformal with respect to the metric $h=g_\tinyZ$ then initializing the flow with $\Phi_\tinyZ$ and evolving for time $t^*+s$ is equivalent to initializing the flow with $\Phi_{t^*}$, setting $h=g_{t^*}$, and evolving  for time $s$. Thus, evolving the limit map under the modified flow must also result in a uniform scaling of the limit map. Since uniform scaling is itself conformal, this implies that the limit map is uniformly scaled by traditional mean-curvature flow and, since the surface is compact, this implies that the surface is a sphere.
\else
($\Leftarrow$) If the map $\Phi_{t^*}\circ\Phi_\tinyZ^{-1}$ is conformal then initializing the flow with surface $\Phi_\tinyZ(M)$ and evolving for time $t^*+s$ is equivalent to initializing the flow with surface $\Phi_{t^*}(M)$ and evolving  for time $s$. Thus, evolving the limit surface under the modified flow must also result in a uniform scaling of the limit surface. Since uniform scaling is itself conformal, this implies that the limit surface is uniformly scaled by traditional mean-curvature flow and, since the surface is compact, this implies that the surface is a sphere.
\fi

\noindent
($\Rightarrow$) If $\Phi_{t^*}$ is a map onto the sphere then we must have
\ifcleanratio
\ifjakeversion
$\Delta_h\Phi_{t^*}=\alpha\sqrt{|h^{\hbox{\tiny -1}}g_{t^*}|}\Phi_{t^*}$
\else
$\Delta_\tinyZ\Phi_{t^*}=\alpha\sqrt{|g_\tinyZ^{\hbox{\tiny -1}}g_{t^*}|}\Phi_{t^*}$
\fi
\else
$\Delta_\tinyZ\Phi_{t^*}=\alpha\sqrt{|g_{t^*}|/|g_\tinyZ|}\Phi_{t^*}$
\fi
\ifjakeversion
for some rescaling constant $\alpha$. In particular, this implies that the heat-flow $\partial\Phi_t/\partial t = \Delta_h \Phi_t$\ will evolve $\Phi_{t^*}$ along directions normal to the surface $\Phi_{t^*}(M)$. Thus the function $\Phi_{t^*}$, considered as a map from $M$ onto the sphere, is harmonic with respect to the metric $h$ and therefore, by Corollary~1 of \cite{Eells:Topology:1976}, conformal.
\else
for some rescaling constant $\alpha$. In particular, this implies that the heat-flow $\partial\Phi_t/\partial t = \Delta_\tinyZ \Phi_t$\ will evolve $\Phi_{t^*}$ along directions normal to the surface $\Phi_{t^*}(M)$. Thus the mapping $\Phi_{t^*}$, considered as a mapping from $M$ to the round sphere, is harmonic with respect to the metric $g_\tinyZ$ and, by Corollary~1 of \cite{Eells:Topology:1976}, therefore conformal.
\fi

\fi

\subsubsection*{What drives it?}
In this work, the flow was derived by modifying traditional mean-curvature flow. However, one can also interpret the flow as a gradient descent on a non-negative energy. In particular, the flow is a descent on the Dirichlet energy of $\Phi_t$ and can be expressed as the sum of a (modified) area energy that drives traditional mean-curvature flow and a conformal energy. (For more details, see Appendix~\ref{a:energy}.)

\subsubsection*{Other geometries}
\ifjakeversion
While the previous discussion considers embeddings of water-tight, genus-zero surfaces, it is also interesting to consider the behavior of the flow on other geometries. To this end, we have applied our flow both to embeddings of surfaces with boundaries, and to embeddings of surfaces with higher genus.
\else
While the previous discussion considers water-tight, genus-zero surfaces, it is also interesting to consider the behavior of the flow on other geometries. To this end, we have applied our flow both to surfaces with boundaries, and surfaces with higher genus.
\fi

\begin{figure}
	\center{\includegraphics[width=\columnwidth]{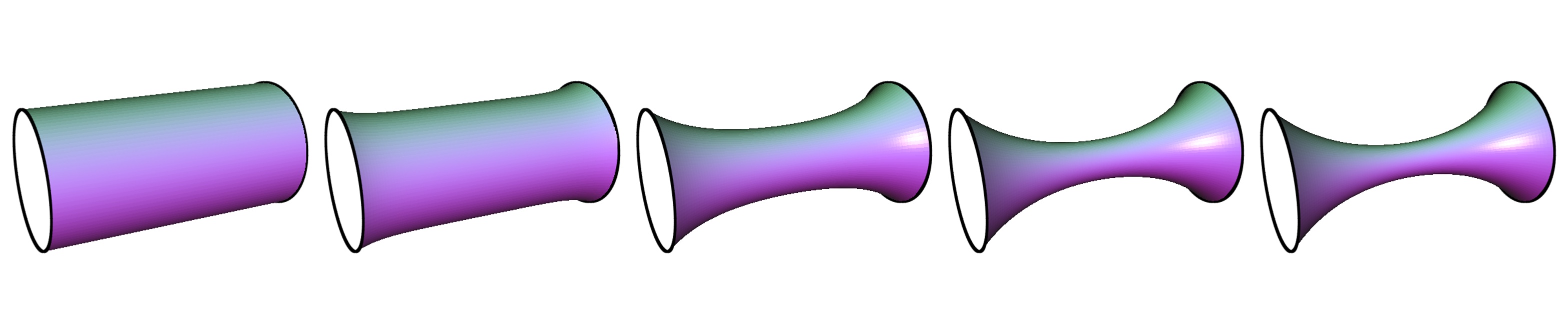}}
  \caption{
  \label{f:cylinder-new}
  Evolution of a cylinder with boundary under the modified mean-curvature flow.
	}
\end{figure}

\begin{figure}
	\center{\includegraphics[width=\columnwidth]{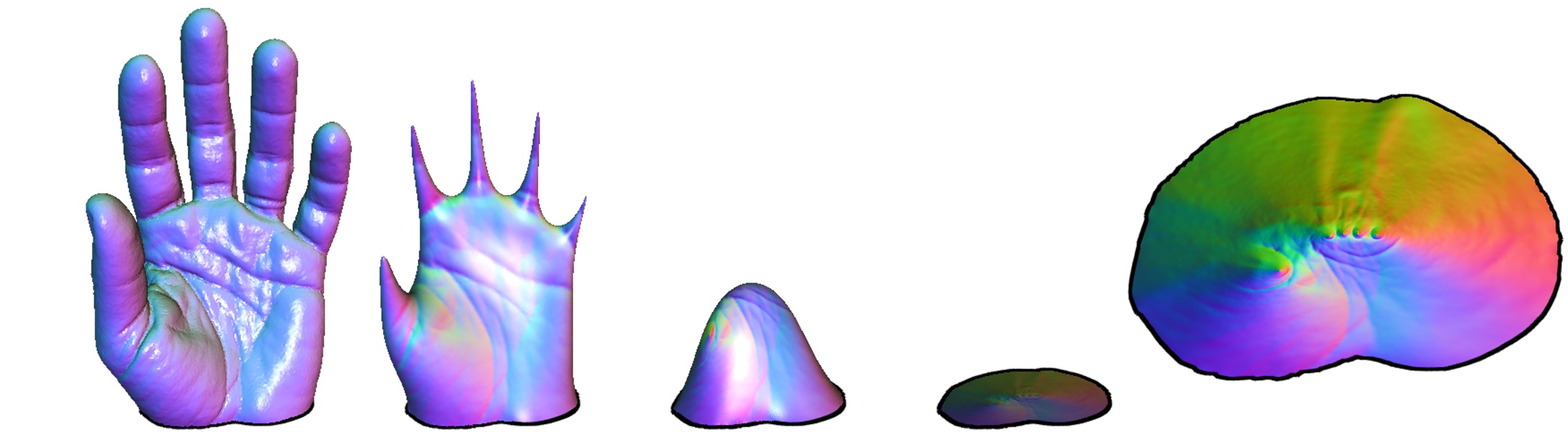}}
  \caption{
  \label{f:hand-new}
  Evolution of a surface with planar boundary under the modified mean-curvature flow (left four images) and a bird's-eye-view of the final mapping (right).
	}
\end{figure}

Figures~\ref{f:cylinder-new} and~\ref{f:hand-new} show the results of our flow for two surfaces with boundaries for which traditional mean-curvature flow generates neck pinches. Though the flow converges, evaluation of the quasi-conformal error shows that the limit map is not conformal. Analyzing the open cylinder, it becomes apparent that stability and conformality cannot be satisfied simultaneously. 
\ifjakeversion
If the mapping were conformal, then the embedding would be fixed under traditional mean-curvature flow, implying that the corresponding surface is minimal. However, since the radius of the cylinder is one and the length of the axis is four, there is no catenoid passing through the two boundaries, and the only minimal surface is the one comprised of two disconnected disks. Thus, the mapping would only be conformal if the flow were to disconnect the surface, which would be impossible without passing through a singularity.
\else
If the mapping were conformal, then the surface would be fixed under traditional mean-curvature flow, implying that it is minimal. However, since the radius of the cylinder is one and the length of the axis is four, there is no catenoid passing through the two boundaries, and the only minimal surface is the one comprised of two disconnected disks. Thus, the mapping would only be conformal if the flow were to disconnect the surface, which would be impossible without passing through a singularity.
\fi


\begin{figure*}
\center
{
	\includegraphics[width=2.1\columnwidth]{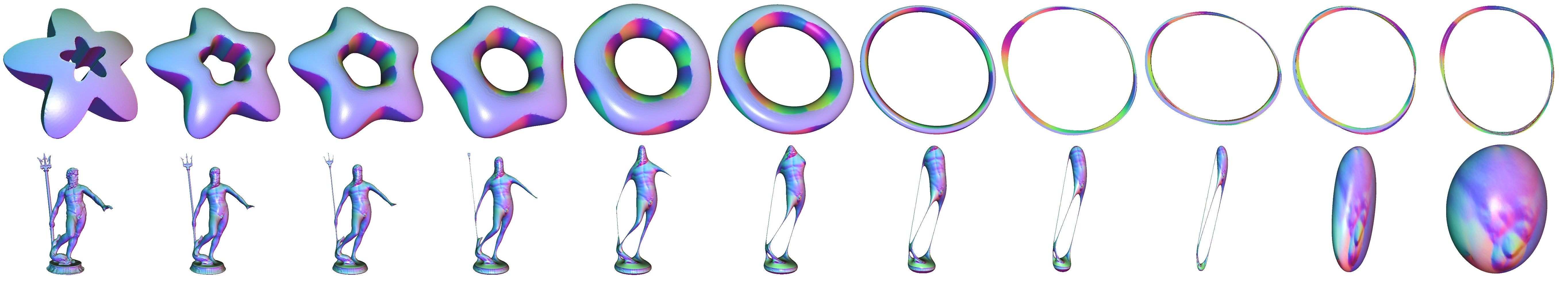}
}
  \caption{
  \label{f:highGenus}
  Visualization of the modified mean-curvature flow on non-genus-zero models.
  }
\end{figure*}

Figure~\ref{f:highGenus} shows the results of our flow on two models that are not simply connected. Though neither flow collapses, they also do not appear to converge. While for higher-genus models the limit map cannot be conformal
\ifjakeversion
(as there are no embeddings of compact, non-genus-zero surfaces that are uniformly scaled by mean-curvature flow) it is possible that a limit map exists. However, even if it does, it is not clear that the limit is a 2-manifold. (For example, the embedding of the star appears to converge to a map onto the circle.) 
\else
(as there are no compact, non-genus-zero surfaces that are uniformly scaled by mean-curvature flow) it is possible that a limit map exists. However, even if it does, it is not clear that the limit is a 2-manifold surface. (For example, the star appears to converge to a circle.) 
\fi

\subsubsection*{Relationship to 1D Flows}
\ifjakeversion
While traditional mean-curvature flow of embeddings of 2D surfaces in 3D can form singularities, this is not the case for embeddings of 1D curves in the plane. In the case of curves, the (uniformly rescaled) flow always converges to a map onto the circle~\cite{Grayson:JDG:1987}. This agrees with the empirical behavior of our modified mean-curvature flow in that the deformation of the 1D curve is always conformal and the definitions of MCF and cMCF agree.

Note that, as in the 2D case, (locally) scaling the map by $\alpha$ scales the Laplace-Beltrami operator by $1/\alpha^2$. However, since the 1D integrals only scale by $\alpha$, the two scaling terms do not cancel out and the discretized Laplace-Beltrami operator does not stay constant.
\else
While traditional mean-curvature flow of 2D surfaces embedded in 3D can form singularities, this is not the case for 1D curves embedded in the plane. In the case of curves, the (uniformly rescaled) flow always converges to a circle~\cite{Grayson:JDG:1987}. This agrees with the empirical behavior of our modified mean-curvature flow in that the deformation of the 1D curve is always conformal and the definitions of the traditional and modified flows agree.

Note that, as in the 2D case, (locally) scaling the geometry by $\alpha$ scales the Laplace-Beltrami operator by $1/\alpha^2$. However, since the 1D integrals only scale by $\alpha$, the two scaling terms do not cancel out and the discretized Laplace-Beltrami operator does not stay constant.
\fi

%% file: conclusion.tex
In this work, we have considered the problem of singularities that arise in mean-curvature flow when evolving non-convex surfaces. Analyzing the finite-elements discretization that commonly arises in geometry processing, we have associated a potential cause for the formation of singularities with the non-conformality of the flow. We have proposed a modification of the flow that simplifies both the discrete and continuous formulations of the flow.  Although we do not have a proof, the work presents empirical evidence that the flow stably evolves genus-zero surfaces, converging to a conformal parameterization of the surface to a sphere.

%% file: appendix_continuous.tex
Although our work has focused on the finite-elements formulation of MCF, the modified flow can also be formulated in a continuous framework.

\ifjakeversion
\noindent\textbf{Claim}: cMCF is driven by the PDE:
\else
\noindent\textbf{Claim}: The modified MCF is driven by the PDE:
\fi
\ifcleanratio
\ifjakeversion
$$\frac{\partial\Phi_t}{\partial t} = \sqrt{\left|\dgh\right|}\Delta_h\Phi_t.$$
\else
$$\frac{\partial\Phi_t}{\partial t} = \sqrt{\left|\dgi\right|}\Delta_\tinyZ\Phi_t.$$
\fi
\fi

\noindent\textbf{Proof}: To show this, we choose a test function $B:M\rightarrow\R$ and consider the Galerkin formulation of the PDE:
\ifcleanratio
\ifjakeversion
$$\mathop{\int}_M\left(\frac{\partial\Phi_t}{\partial t} \cdot B\right)\measure{t} = \mathop{\int}_M\left(\sqrt{\left|\dgh\right|}\Delta_h\Phi_t \cdot B\right)\measure{t}.$$
\else
$$\mathop{\int}_M\left(\frac{\partial\Phi_t}{\partial t} \cdot B\right)\measure{t} = \mathop{\int}_M\left(\sqrt{\left|\dgi\right|}\Delta_\tinyZ\Phi_t \cdot B\right)\measure{t}.$$
\fi
\fi
\ifjakeversion
Note that, because the metric changes over the course of the evolution, the integrands are expressed with respect to the measure $\measure{t}$ not $\measure{h}$.
\else
Note that, because the geometry of $\Phi_t(M)$ changes over the course of the evolution, the Galerkin formulation is expressed with respect to the metric on the surface at time $t$, not the original metric on $M$.
\fi

Applying the change-of-coordinates formula, followed by the Divegence Theorem, we get the weak formulation:
\ifcleanratio
\ifjakeversion
\begin{eqnarray*}
\mathop{\int}_{M}\left(\sqrt{\left|\dhg\right|}\frac{\partial\Phi_t}{\partial t} \cdot B\right)\measure{h} 
&=& \mathop{\int}_{M}\Big(\Delta_h\Phi_t \cdot B\Big)\measure{h}\\
&=& -\mathop{\int}_{M}h\left(\nabla_h\Phi_t , \nabla_h B\right)\measure{h}
\end{eqnarray*}
\else
\begin{eqnarray*}
\mathop{\int}_{M}\left(\sqrt{\left|\dg\right|}\frac{\partial\Phi_t}{\partial t} \cdot B\right)\measure{\tinyZ} 
&=& \mathop{\int}_{M}\Big(\Delta_\tinyZ\Phi_t \cdot B\Big)\measure{\tinyZ}\\
&=& -\mathop{\int}_{M}g_\tinyZ\left(\nabla_\tinyZ\Phi_t , \nabla_\tinyZ B\right)\measure{\tinyZ}
\end{eqnarray*}
\fi
\fi
which gives rise to the discretization with varying mass-matrix (lhs) but fixed Laplace-Beltrami operator (rhs), presented in Section~\ref{s:modifying}.

%% file: appendix_analysis.tex
To better understand the behavior of the flows, we consider three simple surfaces: The (hyperbolic) catenoid, the (elliptic) sphere, and the (parabolic) infinite cylinder.
\ifjakeversion
For each, we analyze the evolution of the embedding of the surface under the actions MCF, heat flow, and cMCF.
\else
For each, we analyze the evolution of the surface under the actions of traditional MCF, heat flow, and the modified flow.
\fi

In our analysis, we use the fact that the Laplacian of the embedding function is the mean-curvature weighted normal, $\Delta\Phi=-H\vec{N}$.

\subsection*{Catenoid}
Since the catenoid has mean-curvature zero, $H=0$, the Laplacian of the embedding is zero and the surface does not evolve under any of the flows.

\subsection*{Sphere}
Due to the rotational symmetry of the sphere, we know that its evolution under all three flows will have the form $\Phi_t(p) = r(t)\cdot \vec{N}(p)$ for some radius function $r(t)$ and normal vector $\vec{N}(p)=p$, so that $\partial\Phi_t/\partial t = r'(t)\cdot\vec{N}(p)$. Without loss of generality, we take the initial radius  to be one.

\paragraph{Traditional MCF}
Using the fact that the mean-curvature of a sphere of radius $r$ is $2/r$, we have:
\begin{eqnarray*}
\frac{\partial\Phi_t}{\partial t} = \Delta_t\Phi_t = -\frac{2\vec{N}}{r}&\Longrightarrow& r'(t)=-\frac{2}{r(t)} \\
&\Longrightarrow& r(t)=\sqrt{1-4t}.
\end{eqnarray*}

\paragraph{Heat Flow}
For the heat flow, the surface is evolved by always using the Laplacian at time $t=0$:
\begin{eqnarray*}
\frac{\partial\Phi_t}{\partial t} = \Delta_\tinyZ\Phi_t = r\cdot\Delta_\tinyZ\Phi_\tinyZ = - 2r\cdot\vec{N}&\Longrightarrow& r'(t)=-2r(t) \\
&\Longrightarrow& r(t)=e^{-2t}.
\end{eqnarray*}

\paragraph{Modified MCF}
\ifjakeversion
Using Equation~\ref{eq:cmcf},
\else
Using Equation~\ref{eq:continuous-modification},
\fi
the evolution of the surface under the modified flow can be described by scaling the heat flow by the reciprocal of the area change:
\begin{eqnarray*}
\frac{\partial\Phi_t}{\partial t} = \frac{-2r\cdot\vec{N}}{r^2}= -\frac{2\vec{N}}{r}&\Longrightarrow& r'(t)=-\frac{2}{r(t)} \\
&\Longrightarrow& r(t)=\sqrt{1-4t}.
\end{eqnarray*}

\subsection*{(Infinite) Cylinder}
Due to the translational and rotational symmetries of the cylinder,
\ifjakeversion
we know that its embedding will evolve
\else
we know that it will evolve
\fi
as a constant offset from the original cylinder along the normal, $\Phi_t(p) = \Phi_\tinyZ(p) - \vec{N}(p) + r(t)\cdot\vec{N}(p)$ for some radius function $r(t)$, so that $\partial\Phi_t/\partial t = r'(t)\cdot\vec{N}(p)$.  Without loss of generality, we take the initial radius to be one.


\paragraph{Traditional MCF}
Since the mean-curvature of a cylinder of radius $r$ is $1/r$, we have:
\begin{eqnarray*}
\frac{\partial\Phi_t}{\partial t} = \Delta_t\Phi_t = -\frac{\vec{N}}{r}&\Longrightarrow& r'(t)=-\frac{1}{r} \\
&\Longrightarrow& r(t)=\sqrt{1-2t}.
\end{eqnarray*}

\paragraph{Heat Flow}
For the heat flow, we use the fact that, on a unit-radius cylinder, the Laplacian of the normal is equal to the Laplacian of the embedding, $\Delta_\tinyZ\vec{N} = \Delta_\tinyZ\Phi_\tinyZ = -\vec{N}$:
\begin{eqnarray*}
\frac{\partial\Phi_t}{\partial t} = \Delta_\tinyZ\left(\Phi_\tinyZ - \vec{N} + r\cdot\vec{N}\right) = - r\cdot\vec{N} &\Longrightarrow& r'(t)=-r(t) \\
&\Longrightarrow& r(t)=e^{-t}.
\end{eqnarray*}

\paragraph{Modified MCF}
\ifjakeversion
As above, we use Equation~\ref{eq:cmcf} to get:
\else
As above, we use Equation~\ref{eq:continuous-modification} to get:
\fi
\begin{eqnarray*}
\frac{\partial\Phi_t}{\partial t} = \frac{-r\cdot\vec{N}}{r}= -\vec{N}&\Longrightarrow& r'(t)=-1 \\
&\Longrightarrow& r(t)=1-t.
\end{eqnarray*}

%% file: appendix_energy.tex
\ifjakeversion
We show that cMCF can be formulated as the gradient flow of an energy that is the sum of a smoothness term, adapted from the energy defining traditional MCF, and a conformal energy. As we will see, this modified energy is just the Dirichlet energy of $\Phi_t$ with respect to the metric $h$.
\else
We show that the modified MCF can be formulated as the gradient flow of an energy that is the sum of a smoothness term, adapted from the energy defining traditional MCF, and a conformality energy. As we will see, this modified energy is just the Dirichlet energy of $\Phi_t$ over $M$.
\fi

\subsection*{Energies}
Using the Euler-Lagrange formulation, traditional MCF can be defined as the gradient flow of the area functional:
\ifcleanratio
\ifjakeversion
$$E_A(\Phi_t) = \hbox{Area}\big(\Phi_t(M)\big)=\mathop{\int}_M\frac{\left|\dhg\right|}{\sqrt{\left|\dhg\right|}}\measure{h}.$$
\else
$$E_A(\Phi_t) = \hbox{Area}\big(\Phi_t(M)\big)=\mathop{\int}_M\frac{\left|\dg\right|}{\sqrt{\left|\dg\right|}}\measure{\tinyZ}.$$
\fi
\fi

\ifjakeversion
We define a conformal energy by measuring the extent to which $\dhg$ differs from a scalar multiple of the identity. Specifically, setting $\hbox{Tr}(\dhg/2)\cdot\hbox{id.}$ to be the scalar multiple of the identity with the same trace as $\dhg$, we set:
\else
We define a conformal energy by measuring the extent to which the ratio of $g_t$ and $g_\tinyZ$ differs from a scalar multiple of the identity. Specifically, setting $\hbox{Tr}(\dg/2)\cdot\hbox{id.}$ to be the scalar multiple of the identity with the same trace as $\dg$, we set:
\fi
\ifcleanratio
\begin{align*}
E_C(\Phi_t)
&=       \mathop{\int}_M\left(\frac{\left\|\dhg - \hbox{Tr}(\dhg/2)\cdot\hbox{id.}\right\|_F^2}{\left|\dhg\right|}\right)\measure{t}\\
&=\frac12\mathop{\int}_M\frac{\hbox{Tr}^2(\dhg) - 4 \left|\dhg\right|}{\sqrt{\left|\dhg\right|}}\measure{h}
\end{align*}
\ifjakeversion
\else
\begin{align*}
E_C(\Phi_t)
&=       \mathop{\int}_M\left(\frac{\left\|\dg - \hbox{Tr}(\dg/2)\cdot\hbox{id.}\right\|_F^2}{\left|\dg\right|}\right)\measure{t}\\
&=\frac12\mathop{\int}_M\frac{\hbox{Tr}^2(\dg) - 4 \left|\dg\right|}{\sqrt{\left|\dg\right|}}\measure{\tinyZ}
\end{align*}
\fi
\fi
\ifjakeversion
where the division by the determinant of $\dhg$ makes the integrand invariant to uniform scaling of the map $\Phi_t$.
\else
where the division by the determinant of $\dg$ makes the integrand invariant to uniform scaling of $\Phi_t(M)$.
\fi


\subsection*{Modifying the Energies}
To define the energy driving our modified MCF, we simplify the energies by replacing the geometric mean of the eigenvalues in the denominator,
\ifcleanratio
\ifjakeversion
$\sqrt{\left|\dhg\right|}$,
\else
$\sqrt{\left|\dg\right|}$,
\fi
\fi
\ifjakeversion
with the arithmetic mean, $\hbox{Tr}(\dhg)/2$,
\else
with the arithmetic mean, $\hbox{Tr}(\dg)/2$,
\fi
and sum the (modified) energies.\footnote{
\ifjakeversion
Note that both denominators scale quadratically with $\Phi_t$, we have $\hbox{Tr}(\dhg/2)\leq\sqrt{\left|\dhg\right|}$, and the two are equal if and only if $\Phi_t$ is conformal with respect to the metric $h$.
\else
Note that both denominators scale quadratically with $\Phi_t$, we have $\hbox{Tr}(\dg/2)\leq\sqrt{\left|\dg\right|}$, and the two are equal if and only if $\Phi_t$ is conformal with respect to the metric $g_\tinyZ$. 
\fi
}

Replacing the denominators, we get:
\ifcleanratio
\ifjakeversion
\begin{eqnarray*}
\tilde{E}_A(\Phi_t) &=& 2\mathop{\int}_M\frac{\left|\dhg\right|}{\hbox{Tr}(\dhg)}\measure{h}\\
\tilde{E}_C(\Phi_t) &=& \frac12\mathop{\int}_M\frac{\hbox{Tr}^2(\dhg) - 4 \left|\dhg\right|}{\hbox{Tr}(\dhg)}\measure{h}.
\end{eqnarray*}
\else
\begin{eqnarray*}
\tilde{E}_A(\Phi_t) &=& 2\mathop{\int}_M\frac{\left|\dg\right|}{\hbox{Tr}(\dg)}\measure{\tinyZ}\\
\tilde{E}_C(\Phi_t) &=& \frac12\mathop{\int}_M\frac{\hbox{Tr}^2(\dg) - 4 \left|\dg\right|}{\hbox{Tr}(\dg)}\measure{\tinyZ}.
\end{eqnarray*}
\fi
\fi
And, taking the sum of the energies, we get:
\ifjakeversion
$$
\tilde{E}(\Phi_t)
= \tilde{E}_A(\Phi_t) + \tilde{E}_C(\Phi_t)
= \frac12\mathop{\int}_Mh(\nabla_h\Phi_t,\nabla_h\Phi_t)\measure{h}.
$$
That is, we replace the area functional defining traditional MCF with the Dirichlet energy of the map $\Phi_t$ with respect to the metric $h$.
\else
$$
\tilde{E}(\Phi_t)
= \tilde{E}_A(\Phi_t) + \tilde{E}_C(\Phi_t)
= \frac12\mathop{\int}_M\left\|\nabla_\tinyZ\Phi_t\right\|_0^2\measure{\tinyZ}.
$$
That is, we replace the area functional defining traditional MCF with the Dirichlet energy of the map $\Phi_t$ with respect to the metric $g_\tinyZ$.
\fi

Linearizing the energy by considering $\Psi:M\rightarrow\R^3$ gives:
\ifcleanratio
\ifjakeversion
\begin{align*}
\tilde{E}(\Phi_t+\epsilon\Psi)
&= \tilde{E}(\Phi_t) - \epsilon\mathop{\int}_M(\Delta_h\Phi_t)\cdot\Psi\measure{h} + O(\epsilon^2)\\
&= \tilde{E}(\Phi_t) - \left\langle\Delta_h\Phi_t,\Psi\right\rangle_h + O(\epsilon^2)\\
&= \tilde{E}(\Phi_t) - \left\langle\sqrt{\left|\dgh\right|}\Delta_h\Phi_t,\Psi\right\rangle_t + O(\epsilon^2),
\end{align*}
\else
\begin{align*}
\tilde{E}(\Phi_t+\epsilon\Psi)
&= \tilde{E}(\Phi_t) - \epsilon\mathop{\int}_M(\Delta_\tinyZ\Phi_t)\cdot\Psi\measure{\tinyZ} + O(\epsilon^2)\\
&= \tilde{E}(\Phi_t) - \left\langle\Delta_\tinyZ\Phi_t,\Psi\right\rangle_\tinyZ + O(\epsilon^2)\\
&= \tilde{E}(\Phi_t) - \left\langle\sqrt{\left|\dgi\right|}\Delta_\tinyZ\Phi_t,\Psi\right\rangle_t + O(\epsilon^2),
\end{align*}
\fi
\fi
\ifjakeversion
where $\langle\cdot,\cdot\rangle_t$ is the inner-product defined on the space of functions on $M$ by the metric $g_t$:
\else
where $\langle\cdot,\cdot\rangle_t$ is the inner-product defined on the space of functions on $M$ by pulling back the metric from $\Phi_t(M)$:
\fi
$$\left\langle f , g\right\rangle_t = \mathop{\int}_{M}(f\cdot g)\measure{t}.$$

\ifcleanratio
\ifjakeversion
Thus, the gradient of the energy $\tilde{E}$, defined with respect to the inner-product $\langle\cdot,\cdot\rangle_t$, is $-\sqrt{\left|\dgh\right|}\Delta_h\Phi_t$.
\else
Thus, the gradient of the energy $\tilde{E}$, defined with respect to the metric on $\Phi_t(M)$, is $-\sqrt{\left|\dgi\right|}\Delta_\tinyZ\Phi_t$.
\fi
\fi

\subsection*{Relationship to Heat Flow}
\ifjakeversion
Note that the heat flow:
$$\frac{\partial\Phi_t}{\partial t}=\Delta_h\Phi_t$$
can also be defined as a gradient flow on the Dirichlet energy of $\Phi_t$ with respect to the metric $h$. However, as demonstrated in Section~\ref{s:results}, heat flow and cMCF evolve the surfaces in different ways. This is because defining the energy gradient requires choosing an inner-product on the space of functions on $M$. For heat flow, the inner-product is defined by the initial metric $h$ giving $\langle\cdot,\cdot\rangle_h$ while for cMCF it is defined by the metric $g_t$ giving $\langle\cdot,\cdot\rangle_t$.

In particular, though both (un-normalized) flows evolve towards the same critical point, (the constant map taking all points in $M$ to a single point in 3D), they converge to this function in different ways. Our modified MCF appears to converge to the constant function ``as a sphere'' (\cite{Huisken:JDG:1984}) while the heat flow does not.
\else
Note that the heat flow:
$$\frac{\partial\Phi_t}{\partial t}=\Delta_\tinyZ\Phi_t$$
can also be defined as a gradient flow on the Dirichlet energy of $\Phi_t$ over $M$. However, as demonstrated in Section~\ref{s:results}, the heat flow and our modified MCF evolve the surfaces in different ways.
This is because the gradient is defined by choosing an area measure on $M$. For the heat flow, the measure is defined by the embedding $\Phi_0(M)$ while for our modified MCF it is defined by the embedding $\Phi_t(M)$.

In particular, though both (un-normalized) flows evolve towards the same critical point, (the constant map taking all points in $M$ to a single point in 3D), they converge to this function in different ways. Our modified MCF appears to converge to the point ``as a sphere'' (\cite{Huisken:JDG:1984}) while the heat flow does not.
\fi

%% file: paper.arxiv.bbl
\begin{thebibliography}{\protect\citename{Springborn et~al\mbox{.} }2008}

\bibitem[\protect\citename{Au et~al\mbox{.} }2008]{Au:SIGGRAPH:2008}
{\sc Au, O., Tai, C., Chu, H., Cohen-Or, D., and Lee, T.}
\newblock 2008.
\newblock Skeleton extraction by mesh contraction.
\newblock {\em ACM Transactions on Graphics (SIGGRAPH '08) 27}, 3.

\bibitem[\protect\citename{Chopp }1993]{Chopp:JCP:1993}
{\sc Chopp, D.}
\newblock 1993.
\newblock Computing minimal surfaces via level set curvature flow.
\newblock {\em Journal of Computational Physics 106\/}, 77--91.

\bibitem[\protect\citename{Crane et~al\mbox{.} }2011]{Crane:SIGGRAPH:2011}
{\sc Crane, K., Pinkall, U., and Schr\"{o}der, P.}
\newblock 2011.
\newblock Spin transformations of discrete surfaces.
\newblock {\em Transactions on Graphics (SIGGRAPH '11) 30\/}.

\bibitem[\protect\citename{Davis and Hager }1999]{David:JMAA:1999}
{\sc Davis, T., and Hager, W.}
\newblock 1999.
\newblock Modifying a sparse {C}holesky factorization.
\newblock {\em SIAM Journal on Matrix Analysis and Applications 20\/},
  606--627.

\bibitem[\protect\citename{Desbrun et~al\mbox{.} }1999]{Desbrun:SIGGRAPH:1999}
{\sc Desbrun, M., Meyer, M., Schr\"{o}der, P., and Barr, A.}
\newblock 1999.
\newblock Implicit fairing of irregular meshes using diffusion and curvature
  flow.
\newblock In {\em ACM SIGGRAPH Conference Proceedings}, 317--324.

\bibitem[\protect\citename{Dziuk }1988]{Dziuk:LNM:1988}
{\sc Dziuk, G.}
\newblock 1988.
\newblock Finite elements for the {B}eltrami operator on arbitrary surfaces.
\newblock In {\em Partial Differential Equations and Calculus of Variations,
  Lecture Notes in Mathematics}, vol.~1357. 142--155.

\bibitem[\protect\citename{Dziuk }1990]{dziuk1990}
{\sc Dziuk, G.}
\newblock 1990.
\newblock An algorithm for evolutionary surfaces.
\newblock {\em Numerische Mathematik 58}, 1, 603--611.

\bibitem[\protect\citename{Eells and Wood }1976]{Eells:Topology:1976}
{\sc Eells, J., and Wood, J.}
\newblock 1976.
\newblock Restrictions on harmonic maps of surfaces.
\newblock {\em Topology 15}, 3, 263--266.

\bibitem[\protect\citename{Grayson }1987]{Grayson:JDG:1987}
{\sc Grayson, M.}
\newblock 1987.
\newblock The heat equation shrinks embedded plane curves to round points.
\newblock {\em Journal of Differential Geometry 26\/}, 285--314.

\bibitem[\protect\citename{Huisken }1984]{Huisken:JDG:1984}
{\sc Huisken, G.}
\newblock 1984.
\newblock Flow by mean curvature of convex surfaces into spheres.
\newblock {\em Journal of Differential Geometry 20\/}, 237--266.

\bibitem[\protect\citename{Mantegazza }2011]{mantegazza}
{\sc Mantegazza, C.}
\newblock 2011.
\newblock {\em Lecture Notes on Mean Curvature Flow}.
\newblock Birkhauser Verlag.

\bibitem[\protect\citename{Pinkall and Polthier }1993]{Pinkall:EM:1993}
{\sc Pinkall, U., and Polthier, K.}
\newblock 1993.
\newblock Computing discrete minimal surfaces and their conjugates.
\newblock {\em Experimental Mathematics 2\/}, 15--36.

\bibitem[\protect\citename{Sander et~al\mbox{.} }2001]{Sander:SIGGRAPH:2001}
{\sc Sander, P., Snyder, J., Gortler, S., and Hoppe, H.}
\newblock 2001.
\newblock Texture mapping progressive meshes.
\newblock In {\em ACM SIGGRAPH Conference Proceedings}, 409--416.

\bibitem[\protect\citename{Springborn et~al\mbox{.}
  }2008]{Springborn:SIGGRAPH:2008}
{\sc Springborn, B., Schr\"{o}der, P., and Pinkall, U.}
\newblock 2008.
\newblock Conformal equivalence of triangle meshes.
\newblock {\em ACM Transactions on Graphics (SIGGRAPH '08) 27\/}, 77:1--77:11.

\bibitem[\protect\citename{Taubin }1995]{Taubin:SIGGRAPH:1995}
{\sc Taubin, G.}
\newblock 1995.
\newblock A signal processing approach to fair surface design.
\newblock In {\em ACM SIGGRAPH Conference Proceedings}, 351--358.

\end{thebibliography}
